\documentclass[11pt,a4paper]{article}
\usepackage{amsmath}
\usepackage{amsfonts}
\usepackage{amssymb}
\usepackage{graphicx}
\usepackage[utf8]{inputenc}
\usepackage{array}
\usepackage{t1enc}
\usepackage[mathscr]{eucal}
\usepackage{amsthm}
\usepackage{mathtools}

\usepackage{color}

\usepackage[colorlinks=false,pdfborderstyle={/W 1}]{hyperref}
 
\def\arXivo#1{\href{http://front.math.ucdavis.edu/#1}{{\tt [arXiv:#1]}}}

\numberwithin{equation}{section}
\newenvironment*{pf}{\noindent{\it Proof.}}{\hfill $\blacksquare$}

\newenvironment*{pfo}{\noindent{\it Proof of}}{\hfill $\blacksquare$}

\theoremstyle{plain}
\newtheorem{tetel}{Theorem}[section]
\theoremstyle{definition}
\newtheorem{defi}[tetel]{Definition}
\theoremstyle{plain}
\newtheorem{prop}[tetel]{Proposition}
\theoremstyle{plain}
\newtheorem{lemma}[tetel]{Lemma}
\theoremstyle{definition}
\newtheorem{rem}[tetel]{Remark}
\theoremstyle{plain}
\newtheorem{cor}[tetel]{Corollary}
\theoremstyle{definition}

\theoremstyle{plain}

\begin{document}

\newcommand{\diam}{\mathop{\textrm{diam}}}
\newcommand{\dist}{{\mathop{\textrm{dist}}}}
\newcommand{\lra}{\leftrightarrow}
\newcommand{\xlra}{\xleftrightarrow}
\newcommand{\xnlra}{\xnleftrightarrow}
\newcommand{\pr}{\mathrm{\mathbb{P}}}
\newcommand{\pp}{\mu}
\newcommand{\ex}{\mathrm{\mathbb{E}}}
\newcommand{\ee}{\mathrm{\overline{\mathbb{E}}}}
\newcommand{\Em}{m_E}
\newcommand{\C}{\mathcal{C}}
\newcommand{\F}{\mathcal{F}}
\newcommand{\A}{\mathcal{A}}
\newcommand{\E}{\mathcal{E}}
\newcommand{\h}{\mathcal{H}}
\newcommand{\om}{{\omega}}
\newcommand{\ebd}{\partial_E}
\newcommand{\ivbd}{\partial_V^\mathrm{in}}
\newcommand{\ovbd}{\partial_V^\mathrm{out}}
\newcommand{\QQ}{\mathcal{Q}}
\newcommand{\w}{{\bf w}}
\newcommand{\dd}{{\bf d}}
\newcommand{\Z}{\mathbb{Z}}
\newcommand{\R}{\mathbb{R}}
\newcommand{\T}{\mathbb{T}}
\newcommand{\N}{\mathbb{N}}
\newcommand{\GG}{\mathcal{G}}
\newcommand{\TT}{\mathcal{T}}
\newcommand{\VV}{\mathcal{V}}
\newcommand{\PP}{\mathcal{P}}
\newcommand{\NN}{\mathcal{N}}
\newcommand{\RR}{\mathcal{R}}
\newcommand{\XX}{\mathcal{X}}
\newcommand{\YY}{\mathcal{Y}}
\newcommand{\B}{\mathcal{B}}
\newcommand{\eps}{\varepsilon}
\newcommand{\cc}{\mathbf{c}}
\newcommand{\degi}{\deg^{in}}
\newcommand{\dego}{\deg^{out}}
\newcommand{\xin}{\underline{\xi}^{(n)}}
\newcommand{\xinpm}{\underline{\xi}^{(n\pm)}}

\title{Controllability, matching ratio and graph convergence}

\author{Dorottya~Beringer\thanks{Alfr\'ed R\'enyi Institute of Mathematics, Hungarian Academy of Sciences, Re\'altanoda u. 13-15, Budapest 1053 Hungary}
\thanks{Supported by the ERC Consolidator Grant 648017.}
 \and
\'Ad\'am~Tim\'ar\footnotemark[1]
\thanks{Supported by the Hungarian National Research, Development and Innovation Office, NKFIH grant K109684, and by grant LP 2016-5 of the Hungarian Academy of Sciences.}
}


\maketitle

\begin{abstract}
There is an important parameter in control theory which is 
closely related to the directed matching ratio of the network, as shown in \cite{LSB}. 
We give proofs on two main statements of the paper of Liu, Slotine and Barab\'asi \cite{LSB} on the directed matching ratio, which were based on numerical results and heuristics from  statistical physics. First, we show that the directed matching ratio of directed random networks given by a fix sequence of degrees is concentrated around its mean. We also examine the convergence of the (directed) matching ratio of a random (directed) graph sequence that converges in the local weak sense, and generalize the result of \cite{EL}. We prove that the mean of the directed matching ratio converges to the properly defined matching ratio parameter of the limiting graph. We further show the almost sure convergence of the matching ratios for the most widely used families of scale-free networks, which was the main motivation of \cite{LSB}.
\end{abstract}

\section{Introduction and results}

Liu, Slotine and Barab\'asi \cite{LSB} examined the controllability of both real networks and network models. 
The models that were most relevant to them are the so-called scale-free networks, which are known to exhibit several characteristics, such as a power-law degree decay, of the networks observed in real-world applications. 
Informally, the controllability parameter of a network is defined as the minimum number $N_D$ of nodes needed to control a network, e.g.
the number of nodes, which can shift molecular networks of the cell from a malignant state to a healthy state. 
They showed that the proportion $n_D=N_D/|V(G)|$ of nodes needed to control a finite network $G$ equals one minus the relative size of the maximal directed matching (directed matching ratio, see Definition \ref{def.dirm}). This allows one to prove results on $n_D$ by proving the corresponding statement for the directed matching ratio. 
In the paper \cite{LSB} it was also observed that the matching ratio is mainly determined by the degree sequence of the graph, namely, if the edges are randomized in a way that does not change the degrees, then the matching ratio does not alter significantly. 
Furthermore, for the most widely used families of scale-free networks, the directed matching ratio converges to a constant. 
These two latter statement were based on numerical results, and for the last one there were also used methods from statistical physics. 
In this paper we give rigorous mathematical proofs of these results on the directed matching ratio. 

Our first theorem gives a quantitative result on the observation that the matching ratio is concentrated if we randomize the edges of a directed graph in a way that does not change the in- and out-degrees. Furthermore, we show that a similar concentration holds if we randomize the edges in such a way that preserves the total degrees but can alter the number of edges pointing to or from the particular vertices. For the definition of the random configuration model used in the next theorem, see Section \ref{ss.networkmodels}. 
\begin{tetel}[Concentration of the matching ratio]\label{thm.concentrationm}
Consider a sequence of in- and out-degrees $d_1^+,\dots ,d_n^+$, respectively $d_1^-,\dots ,d_n^-$, and let $d_j=d_j^++d_j^-$. 

1) Let $G$ be a random directed graph on $n$ vertices given by the random configuration model conditioned on the event that the in- and out-degrees are $d_1^+,\dots ,d_n^+$, respectively $d_1^-,\dots ,d_n^-$. Then the directed matching ratio $m(G)$ of $G$ satisfies 
\begin{align*}
\pr\left(|m(G)-\ex(m(G))|>\eps\right)\leq 2\exp\left\{ -\frac{\eps^2n^2}{8\sum_{k=1}^nd_k^2}. \right\}
\end{align*}

2) Let $G$ be a random directed graph on $n$ vertices given by the random configuration model conditioned on the event that the total degrees of the vertices are $d_1,\dots, d_n$. Then the directed matching ratio $m(G)$ of $G$ satisfies 
\begin{align*}
\pr\left(|m(G)-\ex(m(G))|>\eps\right)\leq 2\exp\left\{ -\frac{\eps^2n^2}{32\sum_{k=1}^nd_k^2}. \right\}
\end{align*}
\end{tetel}

Consider random graph models which ensure a uniform finite bound on the empirical second moments with probability tending to 1. 
Theorem \ref{thm.concentrationm} shows that for graph sequences given by such models, we have a strong concentration of the matching ratio around its mean in the re-randomized graphs with high probability. 
In particular, Erd\H{o}s--R\'enyi graphs or graphs given by the random configuration model with degree distribution $\xi$ with finite second moment have this property. 

Our second result proves the convergence of the matching ratio in the most common families of directed networks. See Definitions \ref{def.lwc} and \ref{def.asconvgr} and Remark \ref{rem.lwc} for the notion of graph convergence and Definition \ref{def.unimodular} for unimodularity. For the graph models used in the theorem see Section \ref{ss.networkmodels}. 
\begin{tetel}[Almost sure convergence of the matching ratio for scale-free graphs]\label{thm.convmall}
1) Let $G_n$ be a sequence of random (directed) finite graphs that converges to a random rooted (directed) graph $(G,o)$ in the local weak sense. 
Then 
\begin{align*}
\lim_{n\to\infty}\ex(m(G_n))=\sup_M\pr_G\big(o\in V^{(-)}(M)\big),
\end{align*}
where the supremum is taken over all (directed) matchings $M$ of $G$ such that the law of $(G,M,o)$ is unimodular.  

2) Let $G_n$ be a sequence of undirected finite graphs defined on a common probability space that converge almost surely in the local weak sense and let $G_n^d$ be a sequence of random directed graphs obtained from $G_n$ by giving each edge a random orientation independently. Then $m(G_n^d)$ converges almost surely to the constant $\lim_{n\to\infty}\ex(m(G_n^d))$. 

3) Let $G_n$ be the sequence of random directed graphs given by the preferential attachment rule. Then $m(G_n)$ converges almost surely to the constant $\lim_{n\to\infty}\ex(m(G_n))$.
\end{tetel}
We prove these results in Section \ref{s.convm}. 
In Subsection \ref{ss.convem} we prove part 1): in Theorem \ref{thm.convem} we show the convergence of the mean of the matching ratio. 
It was proven in \cite{EL} that the limit of the matching ratio of local weak convergent sequences of \emph{deterministic} finite graphs with an uniform bound on the degrees exists. Bordenave, Legrange and Salez \cite{BLS} removed the bounded degree assumption and gave a formula on the value of the limit of the matching ratio. We still need the context of random directed graphs, hence could not apply their result directly. We proceeded through an alternative definition of the matching ratio of the limit object, which looks more natural in our setting. However, the formula in \cite{BLS} for the matching ratio of the limit can be adapted, to obtain quantitative results on the asymptotic value of the directed matching ratio or controllability parameter of large random networks.

In Subsections \ref{sss.diras} we prove the results on the matching ratio that imply part 2). We prove that if a sequence of random directed graphs is obtained from a convergent deterministic graph sequence by orienting each edge independently, then it converges almost surely in the local weak sense, see Definition \ref{def.asconvgr}. 
This is our Lemma \ref{lem.dirasconv} which is similar to Proposition 2.2 in \cite{E09}. 
As a consequence, we get that for directed graphs obtained from almost sure convergent undirected graph sequences the matching ratios converge almost surely. This result applies for sequences given by the random configuration model or Erd\H{o}s--R\'enyi random graphs. 
 
In Subsection \ref{sss.pa} we prove the result that implies part 3) of Theorem \ref{thm.contrconv}. The method used in Subsection \ref{sss.diras} does not apply for the preferential attachment graphs (we cannot start from an a priory almost sure convergence of the undirected graph sequence) hence we needed a different method. 

We note that one can approach the directed matching ratio through an algorithmic point of view, as initiated in \cite{LCsZhP} via the application of the Karp-Sipser algorithm. We do not pursue this direction in the present paper, but preliminary investigations have been started with E. Cs\'oka.

For completeness, we also present our results in the language of controllability. Denote by $n_D$  the proportion of the minimum number of nodes needed to control the network $G$ to the number of nodes, as defined in \cite{LSB}. Our results translate to the following theorems by $n_D(G)=1-m(G)$.

\begin{tetel}[Concentration of the controllability parameter]\label{thm.contrconcentration}
Consider a sequence of in- and out-degrees $d_1^+,\dots ,d_n^+$, respectively $d_1^-,\dots ,d_n^-$, and let $d_j=d_j^++d_j^-$. 

1) Let $G$ be a random directed network on $n$ vertices given by the random configuration model conditioned on the event that the in- and out-degrees are $d_1^+,\dots ,d_n^+$, respectively $d_1^-,\dots ,d_n^-$. Then the controllability parameter $n_D(G)$ of $G$ satisfies 
\begin{align*}
\pr\left(|n_D(G)-\ex(n_D(G))|>\eps\right)\leq 2\exp\left\{ -\frac{\eps^2n^2}{8\sum_{k=1}^nd_k^2}. \right\}
\end{align*}

2) Let $G$ be a random directed network on $n$ vertices given by the random configuration model conditioned on the event that the total degrees of the vertices are $d_1,\dots, d_n$. Then the controllability parameter $n_D(G)$ of $G$ satisfies 
\begin{align*}
\pr\left(|n_D(G)-\ex(n_D(G))|>\eps\right)\leq 2\exp\left\{ -\frac{\eps^2n^2}{32\sum_{k=1}^nd_k^2}. \right\}
\end{align*}
\end{tetel}

\begin{tetel}[Almost sure convergence of the controllability parameter for scale-free graphs]\label{thm.contrconv}
1) Let $G_n$ be a sequence of random directed finite graphs that converges to a random rooted graph $(G,o)$ in the local weak sense. 
Then 
\begin{align*}
\lim_{n\to\infty}\ex(n_D(G_n))=\inf_M\pr_G\big(o\notin V^{(-)}(M)\big),
\end{align*}
where the infimum is taken over all (directed) matchings $M$ of $G$ such that the law of $(G,M,o)$ is unimodular.  

2) Let $G_n$ be a sequence of undirected finite graphs defined on a common probability space that converge almost surely in the local weak sense and let $G_n^d$ be a sequence of random directed graphs obtained from $G_n$ by giving each edge a random orientation independently. Then $n_D(G_n^d)$ converges almost surely to the constant $\lim_{n\to\infty}\ex(n_D(G_n^d))$. 

3) Let $G_n$ be the sequence of random directed graphs given by the preferential attachment rule. Then $n_D(G_n)$ converges almost surely to the constant $\lim_{n\to\infty}\ex(n_D(G_n))$.
\end{tetel}

\subsection{Notations}

We always consider locally finite graphs, with directed or undirected edges. We allow multiple edges and loops. 
We denote by $G\simeq G'$ and $(G,o)\simeq (G',o)$ that the graphs $G$ and $G'$ are isomorphic and rooted isomorphic, respectively. 
We write $\deg_G x$ for the degree of a vertex $x$ in a graph $G$. If the
graph $G$ is directed then denote by $\degi_G x$ and $\dego_G x$ the in- and
out-degree of the vertex $x$. Given a directed edge $e=(x,y)$ we call $x$ the
\emph{tail} and $y$ the \emph{head} of the edge. 
Given a set $F$ of edges let $V(F)$ be the set of vertices that are incident
to an edge in $F$. Let $V^-(F)$, respectively $V^+(F)$ be the set of the tails, respectively
the heads of the edges in $F$. 
Let $B_G(x,n):=\{y\in V(G): \dist_G(x,y)\leq n\}$ be the ball of radius $n$
around a vertex $x$ in the graph $G$ induced by the graph metric. Given a (multi)set
$F$ (of edges or vertices) we denote by $|F|$ the number of elements of the
set (counted with multiplicity). 
Let $[n]$ be the set $\{1,\dots,n\}$. 
Given a random graph $G$ we denote by $\pr_G$ the probability with respect to its law.

\subsection{Directed matchings and graph convergence}\label{ss.dirm}

First we define directed matchings and the matching ratio of directed graphs which are closely related to the controllability of the network. 

\begin{defi}[Directed matching and directed matching ratio]\label{def.dirm}
A \emph{directed matching} $M$ of a directed graph $G$ is a subset of the edges such that the in- and out-degrees in the subgraph induced by $M$ are at most one. 
The \emph{directed matching ratio} of the finite directed graph $G$ is $m(G):=\frac{|V^-(M_{max}(G))|}{|V(G)|}=\frac{|M_{max}(G)|}{|V(G)|}$, where $M_{max}$ is a maximal size directed matching of $G$. 
For undirected finite graphs $G$ we define the matching ratio as $m(G):=\frac{|V(M_{max}(G))|}{|V(G)|}=\frac{2|M_{max}(G)|}{|V(G)|}$, where $M_{max}$ is a maximal size matching of $G$. 
\end{defi}

For possibly disconnected graphs (for instance Erd\H{o}s--R\'enyi graphs or graphs defined by the random configuration model, see Section \ref{ss.networkmodels}), there is another natural way to define the directed matching ratio.
Viewing them as a unimodular random graph, one takes a uniformly chosen random root, and only keeps the \emph{connected component} of this root. Then one could define the matching ratio as the size of the maximal matching of this component divided by the size oft he component.
Contrary to connected graphs, this later definition can give a random variable even if we consider deterministic but disconnected graphs. 
The reason of using Definition \ref{def.dirm} in this paper is coming from our motivating applications in controllability. 
In a finite directed graph the minimum number of nodes needed to control the
network equals the number of vertices that have in-degree 0 in a maximal directed
matching $M_{max}$ (which equals $|V(G)|-|M_{max}(G)|$); see \cite{LSB}. 
We are thus interested in the directed matching ratio
$m(G)$ of a finite directed graph $G$ provided by Definition \ref{def.dirm}, which takes the proportion of vertices of the (possibly disconnected) network that are not needed to control the dynamics of the system. 

In this section we describe the relationship between the matching ratio of directed and undirected graphs. We further define the local weak convergence of graph sequences. 

\begin{defi}[Bipartite representation of a directed graph]\label{def.bipartite}
The \emph{bipartite representation} of a directed graph $G=(V,E)$ is the bipartite graph $G'=(V^-,V^+,E')$ with $V^-=\{v^-: v\in V\}$, $V^+=\{v^+:v\in V\}$ and $E':=\{\{v^-,w^+\}:(v,w)\in E\}$. 
\end{defi}

\begin{rem}\label{rem.dirm}
There is a natural bijection between the directed matchings of $G$ and the matchings of $G'$ which preserves the size of the matching, namely if $M$ is a directed matching of $G$ then $M\mapsto M'=\{\{v^-,w^+\}: (v,w)\in M\}$. Furthermore, $M$ is a directed matching of maximal size if and only if $M'$ is a maximal size matching of $G'$. 
It follows that $m(G)=m(G')$.
\end{rem}

Recall, that a matching $M$ of $G$ has maximal size if and only if there is no augmenting path in $G$ for $M$. By an augmenting path of length $k$ we mean a sequence of disjoint vertices $(v_0,\dots,v_{2k+1})$ such that $\{v_{2j-1},v_{2j}\}\in M$ for $j\in [k]$, $\{v_{2j},v_{2j+1}\}\notin M$ for $j\in \{0,\dots,k\}$ and $\deg_Mv_0=\deg_Mv_{2k+1}=0$.

We examine sequences of networks that have bounded average degrees. 
Benjamini and Schramm \cite{BS} introduced a notion of convergence for such graph sequences: 
\begin{defi}[Local weak convergence of graphs]\label{def.lwc}
We say that the sequence $(G_n,o)$ of locally finite random rooted graphs converge to the locally finite connected random graph $(G,o)$ in the \emph{local weak sense} if for any positive integer $r$ and any finite rooted graph $(H,o)$ we have $\pr\big(B_{G_n}(o,r)\simeq (H,o)\big)\to \pr\big(B_{G}(o,r)\simeq (H,o)\big)$. 
\end{defi}

\begin{rem}\label{rem.lwc}
By the local weak convergence of a sequence $G_n$ of non-rooted finite graphs we always mean the convergence of the sequence with a root chosen uniformly at random among the vertices. 
\end{rem}

For some of the examined graph sequences the following stronger property holds as well: 

\begin{defi}[Almost sure local weak convergence]\label{def.asconvgr} 
Let $G_n$ be a sequence of finite (directed) random graphs defined on a common probability space (if we do not specify the probability space, then we always consider the product space). We say that $G_n$ \emph{converges almost surely in the local weak sense} if almost every realizations of $G_n$ satisfy that the sequence of the deterministic graphs converges in the local weak sense. 
\end{defi}

Finite random graphs with a uniformly chosen root and random rooted graphs that are local weak limits of (random) finite graphs, satisfy the so-called Mass Transport Principle, see \cite{BS}, Section 3.2. The class of graphs that obeys this principle are called \emph{unimodular} graphs. 
\begin{defi}[Unimodular graphs]\label{def.unimodular}
A random rooted (directed, labeled) graph $(G,o)$ is called {unimodular} if it obeys the Mass Transport Principle: for every measurable real valued function $f$ on the class of locally finite graphs with an ordered pair of vertices that satisfies $f(G,x,y)=f(\gamma G, \gamma x, \gamma y)$ for every $\gamma\in Aut(G)$ the following holds: 
\begin{align*}
\ex\left(\sum_{x\sim o}f(G,o,x)\right)=\ex\left(\sum_{x\sim o}f(G,x,o)\right).
\end{align*}
\end{defi}

Directed matchings and hence the matching ratio of a finite directed graph $G$ can be examined using the bipartite representation $G'$ as mentioned in Remark \ref{rem.dirm}. 
In the next proposition, we analyze the relationship between a convergent graph sequence and its bipartite representation. 

\begin{prop}\label{prop.bipartite}
If a sequence $G_n$ of random directed graphs converges to the random rooted directed graph $(G,o)$, then 
the bipartite representations $G'_n$ converge to $(G',o')$, where $G'$ is the bipartite representation of $G$ with root $o'$ being $o^-$ or $o^+$ with probability 1/2-1/2. 

The converse does not hold: the convergence of the sequence of bipartite representations $G'_n$ does not imply the convergence of $G_n$.  
In fact, there are different random directed rooted graphs $(G_1,o_1)$ and $(G_2,o_2)$ that are limits of sequences of finite random rooted graphs such that $(G'_1,o'_1)$ is isomorphic to $(G'_2,o'_2)$.
\end{prop}

\begin{pf}
Denote by $\mu_{n,r}$ and $\mu_r$ the law of $B_{G_n}(o,r)$, respectively $B_{G}(o,r)$ in the space of locally finite rooted directed graphs and 
let $\mu'_{n,r}$ and $\mu'_r$ the law of $B_{G'_n}(o',r)$, respectively $B_{G'}(o',r)$ in the space of locally finite rooted graphs. 
The random uniform root $o'$ of a bipartite representation $G'_n$ of a finite directed graph $G_n$ is $o^-$ or $o^+$ with probability 1/2-1/2, where $o$ is a uniform random root of $G$. 
It follows that $\mu'_{n,r}=1/2\mu'_{n,r,o^-}+1/2\mu'_{n,r,o^+}$, where $\mu'_{n,r,o^-}$ and $\mu'_{n,r,o^+}$ are the laws of $B_{G'_n}(o^-,r)$, respectively $B_{G'_n}(o^+,r)$. The first statement of the remark follows. 

An example to the second statement is the following. Let $G_1$ be the graph with vertex set $V(G_1)=\Z$ and edge set $E(G_1)=\{(2k,2k-1), (2k,2k+1):k\in\Z\}$, i.e. the usual graph of $\Z$ with an alternating orientation to the edges. Let the random root $o_1$ be $2k$ or $2l-1$ for some $k,l\in\Z$ with probability 1/2-1/2 (the isomorphism class of $(G_1,o)$ does not depend on the actual choice of the integers $k$ and $l$). This graph is the limit of the cycles $C_{2n}$ with $2n$ vertices and edges with alternating orientations. 
Let $G_2$ be the one-point graph without edges with probability 1/2 and with probability 1/2 let $G_2$ be the infinite regular tree with in- and out-degrees 2. This graph is the limit of the sequence of random graphs on $n$ vertices where with probability 1/2 there are no edges and with probability 1/2 the graph is uniformly randomly chosen from the set of graphs on $n$ vertices with all in- and out-degrees 2. 
Then $(G'_1,o'_1)$ and $(G'_2,o'_2)$ are both isomorphic to the random graph that is the one-point graph without edges or $\Z$ with probability 1/2-1/2. 
\end{pf}

\subsection{Canonical network models and their limits}\label{ss.networkmodels}

Some of the examined graph sequences converge to the so-called unimodular Galton--Watson tree. 

\begin{defi}[Unimodular Galton--Watson tree]\label{def.ugw}
Let $\xi$ be a non-negative integer valued random variable with $\ex\xi<\infty$. The unimodular Galton--Watson tree with offspring distribution $\xi$ (denoted by $UGW(\xi)$) is a random rooted tree with root $o$. 
We say that a vertex $y$ is the child of the vertex $x$, if they are adjacent and $\dist(y,o)=\dist(x,o)+1$.
The graph $UGW(\xi)$ is given by the following recursive definition:  
\begin{itemize}
\item The probability that $o$ has $k\geq 0$ children is $\pr(\xi=k)$.
\item For each vertex $x$ the probability that $x$ has $k\geq 0$ children is $\frac{(k+1)\pr(\xi=k+1)}{\ex\xi}$.
\end{itemize}

Let the directed unimodular Galton--Watson tree $UGW^d(\xi)$ be the random rooted directed graph obtained from $UGW(\xi)$ by orienting each edge independently. 
\end{defi}

Now we present the network models examined in this paper. 
For each model first we define the non-directed model and present the known results on the local weak limit of the sequence, then we give the definition of the directed versions and the local weak limit of them. 

\medskip
\noindent\textbf{Random $d$-regular graphs}

Let $G_n$ be the random graph chosen uniformly at random from the set of graphs on the vertex set $[n]$ with all degrees equal $d$. 
It is standard, that the local weak limit of $G_n$ as $n\to\infty$ is the infinite $d$-regular tree $\T_d$. 
In fact, the random graphs $G_n$ converge almost surely to $\T_d$. This follows from the almost sure convergence of the more general class of graphs given by the random configuration model. 

There are two natural ways to define random directed regular graphs. 
The first one is if $G_n$ is a uniformly chosen directed graph on $[n]$ such that each vertex has in- and out-degrees $d$. The local weak limit is a regular tree with in- and out-degrees $d$. 
The second way to define directed graphs $G_n$ is if we choose a uniform random non-directed $d$-regular graph on $[n]$ and orient each edge uniformly at random independently from each other. This model is a special case of the random configuration model defined in the sequel. 
The limit of that graph sequence is the $d$-regular tree with independently oriented edges. 

\medskip
\noindent\textbf{Erd\H{o}s--R\'enyi random graphs}

The Erd\H{o}s--R\'enyi random graphs $\mathcal{G}_{n,c/n}$ are defined in the following way: consider the complete graph on $n$ vertices and keep each edge with probability $c/n$, and delete each edge with probability $1-c/n$ independently from each other. The resulting random graph is $\mathcal{G}_{n,c/n}$. 

The local weak limit of $\mathcal{G}_{n,c/n}$ is $UGW($Poisson$(c))$, that is the Galton--Watson tree with Poisson($c$) offspring distribution. In fact, for almost every realization of the sequence $\mathcal{G}_{n,c/n}$, that sequence of deterministic graphs converges to $UGW($Poisson$(c))$ as well, see Theorem 3.23 in \cite{Brg}.

We define the directed Erd\H{o}s--R\'enyi random graphs $\mathcal{G}_{n,c/n}^d$ by orienting each edge of $\mathcal{G}_{n,c/n}$ uniformly at random independently for the edges. The local weak limit of this sequence is $UGW^d($Poisson$(c))$. 

\medskip
The next two graphs have become increasingly important in applications, because they grab important characteristics of real-world networks (scale-free networks). This is the reason why in \cite{LSB}, which was motivated by applications of controllability, these graphs were studied.

\medskip
\noindent\textbf{Random configuration model}

We fix a non-negative integer valued probability distribution $\xi$. 
We define the graph $G_n$ in the following way:
let $\xi_1,\dots ,\xi_n$ be i.i.d. variables with distribution $\xi$. 
Given $\xi_1,\dots ,\xi_n$ let $\E:=\{(k,j):k\in[n], j\in[\xi_k]\}$ be the set of the \emph{half-edges}. 
Let $H$ be a uniform random perfect matching of the set $\E$ (if $|\E|$ is odd, then put off one half-edge uniformly at random before choosing a perfect matching). Then $H$ defines the random graph $G_n=G_n(H)$ on $[n]$. 

If $\ex(\xi^2)<\infty$, then $G_n$ converge to $UGW(\xi)$ in the local weak sense (see Theorem 3.15 in \cite{Brg}). Furthermore, if $\ex(\xi^p)<\infty$ with some $p>2$, then for almost every realization of the sequence $G_n$, the local weak limit of that deterministic graph sequence is $UGW(\xi)$; see Theorem 3.28 in \cite{Brg} and Theorem \ref{thm.Bordenave}. 

If we want to define a directed graph, then we orient each edge uniformly at
random independently from the other edges. 
We get the same distribution if after fixing the degree sequence $\xi_1,\dots ,\xi_n$ 
we select a subset $\E_T\subseteq\E$ of size $\lfloor|\E|/2\rfloor$ uniformly at random. 
Then we set $\xi_k^-:=|\{j\in[\xi_k]: (k,j)\in\E_T\}|$, $\xi_k^+:=\xi_k-\xi_k^-$ and 
we denote by $\TT:=\{(k,j,-):k\in[n], j\in[\xi_k^-]\}$ the set of the tail-type half-edges and 
by $\h:=\{(k,j,+):k\in[n], j\in[\xi_k^+]\}$ the set of the head-type half-edges. 
Let $\NN$ be the set of the perfect matchings of $\TT$ to $\h$ and 
denote by $N$ a uniform random element of $\NN$. 
Then $N$ defines the random directed graph $G_n=G_n(N)$ on the vertex set $[n]$.

\medskip
\noindent\textbf{Preferential attachment graphs}

The notion of preferential attachment graphs was introduced by Barab\'asi and
Albert in \cite{BA} and the precise construction was given by  Bollob\'as and
Riordan in \cite{BR}. There are several versions of the definition of this
family of random graphs which have turned out to be asymptotically the same: they all converge to the same infinite limit graph; see \cite{BBCS}. 
Altough in the original definitions the preferential attachment graphs are not directed, 
there is a natural way to give each edge an orientation and these orientations extend to the limit graph as well. 

We will use the following definition from \cite{BBCS} completed with the natural orientation of the edges: 
fix a positive integer $r$ and $\alpha\in[0,1)$. For each $n$ the random graph $G_n=G_{r,\alpha,n}^{P\!A}$ is a graph on the vertex set $[n]$ defined by the following recursion: 
let $G_{0}$ be the graph with one vertex and no edges. Given $G_{n-1}$ we construct $G_{n}$ by adding the new vertex $n$ and $r$ new edges with tails $n$. We choose the heads $w_1,\dots w_r$ of the new edges independently from each other in the following way: with probability $\alpha$ we choose $w_j$ uniformly at random among $[n-1]$, and with probability $1-\alpha$ we choose $w_j$ proportional to $\deg_{G_{n-1}}$. Note that each vertex except the starting vertex has out-degree $r$ and each vertex has a random in-degree with mean converging to $r$. 

Berger, Borgs, Chayes and Saberi proved in \cite{BBCS} that the local weak
limit of $G_{r,\alpha,n}^{P\!A}$ as $n\to \infty$ is the P\'olya-point graph with parameters $r$ and $\alpha$. This graph is a unimodular random infinite tree with directed edges; see \cite{BBCS}, Section 2.3 for the definition.

\section{Concentration of the matching ratio in randomized networks}\label{ss.concentration}

In this section we prove Theorem \ref{thm.concentrationm}, which gives a quantitative version of the following experimental observation of Liu, Slotine and Barab\'asi in \cite{LSB}: if we consider a large directed graph, and randomize the edges in such a way that does not change the in- and out-degrees of the graph, then the matching ratio does not alter significantly. 
Part 1) of Theorem \ref{thm.concentrationm} shows the concentration for randomized graphs with the in- and out-degrees left unchanged. This is the result that was observed through simulations in \cite{LSB}. 
Part 2) of the theorem shows that a very similar concentration phenomenon holds even after a randomizing that does not require the in- and out-degrees to be unchanged but only the total degree to remain the same for every vertex. 
In particular, Theorem \ref{thm.concentrationm} shows that if a graph sequence satisfies that the empirical second moment of the degree sequence is $o(n)$ with probability tending to 1 (as $n\to\infty$), then the directed matching ratios of the graphs with randomized edges are concentrated around their mean. 

First we need a lemma that shows that modifying a (directed) graph just around a few vertices cannot alter the size of the maximal matching too much. 

\begin{lemma}\label{lem.adde}
Adding some new edges with a common endpoint to an undirected finite graph or 
adding edges with a common head (respectively tail) to a directed finite graph can increase the size of the maximal
matching by at most one. 
\end{lemma}

\begin{pf} 
For directed graphs the statement follows from the undirected case, using the bipartite representation (see Definition \ref{def.bipartite}). 
For undirected graphs let $F$ be the set of new edges with common endpoint $x$ and let $G_2$ be the graph with vertex set $V(G)$ and edge set $E(G_2)=E(G)\cup F$. If $M_2$ is a maximal size directed matching of $G_2$, then there is at most one edge in $M_2\cap F$ by the definition of the matching. Then $M_2\setminus F$ is a matching of $G$, hence $|M_{max}(G)|\geq |M_2|-1$. 
\end{pf}

Before proving the proposition, we state a version of the Azuma--Hoeffding
inequality (see \cite{LP}, Theorem 13.2), that we will use in this paper. 

\begin{tetel}[Azuma--Hoeffding inequality]\label{thm.AH}
Let $X_1,\dots ,X_n$ be a series of martingale differences. Then 
\begin{align*}
\pr\left(\sum_{k=1}^n X_k>\eps\right)\leq \frac{\eps^2}{2\sum_{k=1}^n \|X_k\|^2_{\infty}}.
\end{align*}
\end{tetel}

The proof of Theorem \ref{thm.concentrationm} uses similar methods to that of  Corollary 3.27 in \cite{Brg}, which implies the concentration of matching ratio for undirected graphs. 

\medskip
\begin{pfo}\emph{ Theorem \ref{thm.concentrationm}.}
We prove both parts of the theorem in the following way: we define random variables $X_k$, $k\in [n]$ which form a series of martingale differences and satisfy $\sum_{k=1}^nX_k=n(m(G)-\ex(m(G)))$. 
We will show that there is an almost sure bound $|X_k|\leq cd_k$, hence we have by the Azuma--Hoeffding inequality 
\begin{align*}
\pr\left( |m(G(N))-\ex(m(G(N)))|>\eps \right)& 
=\pr\left( |X_1+\dots +X_n|>\eps n \right)\\
&\leq 2\exp\left\{ -\frac{(\eps n)^2}{2\sum_{k=1}^n \|X_k\|_\infty^2} \right\}\\
&\leq 2\exp\left\{ -\frac{\eps^2n^2}{2c^2\sum_{k=1}^{n}d_k^2} \right\}.
\end{align*}

\noindent\emph{Part 1).} 
Recall the second definition of the directed random configuration model from Section \ref{ss.networkmodels}, conditioned on the fixed sequences of in- and out-degrees. 
For a half-edge $h=(i,j,\pm)\in \TT\cup\h$ let $v(h):=i$ be the corresponding
vertex and let $N(h)$ be the pair of the half-edge $h$ by the matching $N$. 
Denote by $N(k):=\{(h,h')\in N: v(h),v(h')\in [k]\}$ the partial matching that
consists of the pairs of half-edges of $N$ with corresponding vertices both in $[k]$. 
Let 
\begin{align}
X_k:=\ex\left(|M_{max}(G(N))|\Big|N(k)\right)-\ex\left(|M_{max}(G(N))|\Big|N(k-1)\right). 
\end{align}
The variables $X_k$ clearly form a series of martingale differences, and we claim that $|X_k|\leq 2d_k$ almost surely for all $k\in[n]$. 

We will show that if
$N_1$ and $N_2$ are two partial matchings of $\TT(k):=\{(l,j,-):l\in[k],j\in [d_l^-]\}$ to $\h(k):=\{(l,j,+):l\in[k],j\in [d_l^+]\}$ 
such that they only differ by an edge with tail $k$, i.e. $N_2=N_1\cup{e}$ with $v(e^-)=k$, then 
\begin{align}\label{eq.G12}
\left|\ex\left(|M_{max}(G(N))|\Big|N(k)=N_1\right)-\ex\left(|M_{max}(G(N))|\Big|N(k)=N_2\right)\right|\leq 2,
\end{align}
and the same holds if $N_1$ and $N_2$ differ only by an edge with head $k$. 
It follows that for any two partial matchings $N_1$ and $N_2$ of $\TT(k)$ to $\h(k)$ that satisfy $N_1(k-1)=N_2(k-1)$ 
the left hand side of \eqref{eq.G12} is at most $4d_k$. This implies the bound on $X_k$.
 
To show \eqref{eq.G12}, 
we fix two arbitrary partial matchings $N_1$ and $N_2$ of $\TT(k)$ to $\h(k)$ such that $N_1(k-1)=N_2(k-1)$ and
$N_2=N_1\cup\{(h,h')\}$ with $v(h)=k$. 
Let $\NN_i:=\{N:N(k)=N_i\}$ for $i=1,2$ be the set of perfect matchings of $\h$ to $\TT$ with
$N(k)=N_i$. 
For a configuration $N\in\NN_1$ let 
\begin{align}
f(N):=\left(N\setminus
\{(h,N(h)),(N(h'),h')\}\right)\cup\{(h,h'),(N(h'),N(h))\}.
\end{align}
For each $N\in\NN_1$ there is a unique $f(N)\in\NN_2$ and for all $N'\in
\NN_2$ the size of the set $\{N\in\NN_1: f(N)=N'\}$ is equal, namely 
$\left(\sum_{j=k+1}^{n}d_j^-\right)-\left(\sum_{j=1}^{k}d_j^+-|N_2|\right)=
\frac{|\NN_1|}{|\NN_2|}$.
We have 
\begin{align}
\Big|\ex&\Big(|M_{max}(G(N))|\big|N(k)=N_1\Big)- \ex\Big(|M_{max}(G(N))|\big|N(k)=N_2\Big) \Big| \nonumber \\
&\leq \sum_{H\in\NN_2}\Big| \ex\left(|M_{max}(G(N))|\big| N\in\NN_1,f(N)=H\right)\pr\left(f(N)=H\big|N\in\NN_1\right)- \nonumber \\
&\hskip 1.2 in \ex\left(|M_{max}(G(N))|\big| N\in\NN_2, N=H\right)\pr\left(N=H\big|N\in\NN_2\right) \Big| \nonumber \\
&= \sum_{H\in\NN_2} 
\left| \ex\left(|M_{max}(G(N))|\big|N\in\NN_1, f(N)=H\right)-|M_{max}(G(H))| \right|
\frac{1}{|\NN_2|}. \label{eq.N}
\end{align}
For any $N\in\NN_1$ with $f(N)=H$ the graphs $G(N)$ and $G(H)$ differ by at most four edges in such a way that the size of the set of the heads of these vertices is at most two. By Lemma \ref{lem.adde} we have in this case 
\begin{align*}
\left| \ex\left(|M_{max}(G(N))|\big|N(k)=N_1, f(N)=H\right)-|M_{max}(G(H))| \right|\leq 2
\end{align*}
which combined with \eqref{eq.N} proves inequality \eqref{eq.G12}. 

\medskip
\noindent\emph{Part 2).}
Recall the notations and the second definition of the directed random configuration model from Section \ref{ss.networkmodels}, conditioned on the fixed sequence of total degrees. 
Let $\E(k):=\{(j,l)\in \E: j\in[k]\}$ consist of all half-edges whose end-vertex is in $[k]$, and similarly for any subset $H\subseteq \E$ let $H(k):=\{(j,l)\in H: j\in[k]\}$. 
We claim that for any fixed $k$ and $j\in [d_k]$, if $F_1$ and $F_2$ are subsets of $\E(k)$ such that $F_2=F_1\cup\{(k,j)\}$, then 
\begin{align}\label{eq.F1F2}
\Bigg| \ex\left(|M_{max}(G_n)|\Big|\E_T(k)=F_1\right)-
\ex\left(|M_{max}(G_n)|\Big|\E_T(k)=F_2\right) \Bigg|\leq 4. 
\end{align}
Let $\F_i:=\{H_i\subseteq\E: |H_i|=|\E|/2, H_i(k)=F_i\}$ for $i=1,2$ and let 
\begin{align*}
\RR:=\{(H_1,H_2)\in\F_1\times \F_2: |H_1\bigtriangleup H_2|=2\}. 
\end{align*}
For every $H_1\in\F_1$ the size of the set $\{H_2: (H_1,H_2)\in\RR\}$ equals $|\RR|/|\F_1|=|\E|/2-|F_1|$ and 
for every $H_2\in\F_2$ the size of the set $\{H_1: (H_1,H_2)\in\RR\}$ equals $|\RR|/|\F_2|=\sum_{j=k+1}^nd_j-\left(|\E|/2-|F_2|\right)$. 
The left hand side of \eqref{eq.F1F2} can be bounded above by 
\begin{align*}
\frac{1}{|\RR|} \sum_{(H_1,H_2)\in \RR} 
\Bigg| \ex\left(|M_{max}(G_n)|\Big|\E_T=H_1\right)-
\ex\left(|M_{max}(G_n)|\Big|\E_T=H_2\right) \Bigg|, 
\end{align*}
where each term in the sum is bounded above by 4 by the following argument. 
Fix $(H_1,H_2)\in\RR$, let $\TT_i$ and $\h_i$ be the set of tail- and head-type
half-edges given by $\E_T=H_i$ for $i=1,2$. 
Let $h_1:=\h_1\setminus \h_2$, $h_2:=\h_2\setminus \h_1$, 
$t_1:=\TT_1\setminus \TT_2$ and $t_2:=\TT_2\setminus \TT_1$. 
For each perfect matching $N\in\NN_1$, let 
\begin{align*}
f(N):=\left(N\setminus\Big\{\Big(t_1,N(t_1)\Big),\Big(N(h_1),h_1\Big)\Big\}\right)
\cup\Big\{\Big(t_2,N(t_1)\Big),\Big(N(h_1),h_2\Big)\Big\},
\end{align*}
which is an element of $\NN_2$. Note that $f:\NN_1\to\NN_2$ is a bijection and 
$G(N)$ and $G(f(N))$ differ by at most 4 edges, hence by Lemma \ref{lem.adde} the
size of the maximum matchings of them differ by at most 4. It follows that 
\begin{align*}
\Bigg| \ex\left(|M_{max}(G_n)|\Big|\E_T=H_1\right)-
\ex\left(|M_{max}(G_n)|\Big|\E_T=H_2\right) \Bigg|&\leq \\
\sum_{N\in\NN_1}\frac{1}{|\NN_1|}  
\bigg| \big|M_{max}(G_n(N))\big|-
\big|M_{max}\big(G_n(f(N))\big)\big| \bigg|
&\leq 4.
\end{align*}
This proves \eqref{eq.F1F2}. 

Let 
\begin{align}
X_k:=\ex\left(|M_{max}(G_n)|\Big|\E_T(k)\right)-
\ex\left(|M_{max}(G_n)|\Big|\E_T(k-1)\right).  
\end{align}
We claim that $|X_k|\leq 4d_k$ almost surely for all $k\in[n]$. 
For any $F\subseteq \E(k)$, let $r(F):=\{(j,l): j\in[k], l\leq |\{i:(j,i)\in F\}|\}$,  
i.e. we transform $F$ to a subset with the same size but with the
smallest possible second coordinates. This transform does not change the isomorphism class of the induced directed graph, hence  
$\ex\left(|M_{max}(G_n)|\Big|\E_T(k)=F\right)
=\ex\left(|M_{max}(G_n)|\Big|\E_T(k)=r(F)\right)$.  
This implies that for any two subsets $F_1$ and $F_2$ of $\E(k)$ with
$F_1(k-1)=F_2(k-1)$, the subsets $r(F_1)$ and $r(F_2)$ differ by at most $d_k$
half-edges that all have first coordinate $k$. It follows by \eqref{eq.F1F2} that 
\begin{align*}
\Bigg| \ex\left(|M_{max}(G_n)|\Big|\E_T(k)=F_1\right)-
\ex\left(|M_{max}(G_n)|\Big|\E_T(k)=F_2\right) \Bigg| &=\\
\Bigg| \ex\left(|M_{max}(G_n)|\Big|\E_T(k)=r(F_1)\right)-
\ex\left(|M_{max}(G_n)|\Big|\E_T(k)=r(F_2)\right) \Bigg|&\leq 4d_k, 
\end{align*}
which implies $|X_k|\leq 4d_k$.
\end{pfo}

\section{Convergence of the matching ratio}\label{s.convm}

The goal of this section is to prove the convergence of the directed matching ratio for convergent sequences of random directed graphs. 
This convergence is understood in the stronger sense of almost sure convergence, as we will see, but the proof will often proceed through showing convergence in expectation and then concentration.
For a fixed deterministic non-directed graph sequence that is locally convergent when a uniform root is taken, the convergence of the matching ratio is proved by Elek and Lippner in \cite{EL} if there is uniform bound on the degrees and by Bordenave, Lelarge and Salez in \cite{BLS} in the unbounded case. 
To prove the results of Liu, Slotine and Barabási in \cite{LSB}, we need to generalize these results for directed random graphs. 

In Subsection \ref{ss.convem} we use the method of Elek and Lippner to prove Theorem \ref{thm.convem} on the convergence of the \emph{expected value} of the directed matching ratio of sequences of random graphs. In Definition \ref{def.infinitem} we give an extension of the definition of the expected matching ratio to unimodular random rooted graphs. 
By Theorem 1 in \cite{BLS} and our Theorem \ref{thm.convem} our definition of the expected matching ratio equals twice the parameter $\gamma$ defined in \cite{BLS}. 

In Subsection \ref{ss.asconvm} we prove the almost sure convergence of the directed matching ratios for the network models defined in Subsection \ref{ss.networkmodels}.

\subsection{Convergence of the mean of the matching ratio}\label{ss.convem}

Elek and Lippner proved that the non-directed matching ratio converges if $G_n$ is a convergent sequence of finite deterministic graphs with uniformly bounded degree; see \cite{EL}, Theorem 1.
There are three properties of our examined models, that do not let us apply
this theorem directly: our graphs do not have bounded degrees, and they
are directed and random graphs. 
Although the degrees are not bounded in the examined models of convergent graph sequences, the expected value of the degree of the uniform random root of the random graphs has a uniform bound in each model. 
In Theorem \ref{thm.convem} we prove the convergence of the mean of the matching ratio for convergent sequences of random directed graphs using the method of Elek and Lippner. 

One can extend the (expected) matching ratio to the class of unimodular random (directed) graphs in a natural way. For finite random graphs, the following definition gives the expected value of the matching ratio. 

\begin{defi}[Matching ratio of an infinite graph and unimodular matchings]\label{def.infinitem}
Let $(G,o)$ be a unimodular random (directed) rooted graph. Then the \emph{(expected) matching ratio} of $(G,o)$ is 
\[\Em(G,o)=\sup_{M}\pr_G(o\in V^{(-)}(M)),\]
where the supremum is taken over all random (directed) matchings of $G$ such that the law of $(G,{M},o)$ is unimodular. Matchings with this property will be called \emph{unimodular matchings}.
\end{defi}

\begin{rem}\label{rem.infinitem}
Let $(G,o)$ be a random directed rooted unimodular graph and let $(G',o')$ be its bipartite representation (see Definition \ref{def.bipartite}). Then Lemma \ref{lem.unimodsup} will imply that $\Em(G,o)=\Em(G',o')$. 
\end{rem}

\begin{tetel}\label{thm.convem} 
Let $G_n$ be a sequence of random finite (directed) graphs that converges to the random (directed) rooted graph $(G,o)$ that has finite expected degree. Then 
\begin{align*}
\lim_{n\to\infty}\ex(m(G_n))=\Em(G,o).
\end{align*}
\end{tetel}

To prove Theorem \ref{thm.convem}, we follow the method of \cite{EL}. 
The main differences to that proof come from the lack of uniform bound on the degrees. 
We will define the matchings $M(T)$ in Lemma \ref{lem.fIIDmatching} as factor of IID, which helps us handle the case of unbounded degrees. 
For graphs with unbounded degrees, Lemma 4.1 of \cite{EL} does not apply, hence we will have to proceed through Lemma \ref{lem.muH}. 

\begin{defi}[Factor of IID]\label{def.fIID}
Let $\GG_\star$ be the set of the isomorphism classes of locally finite rooted (directed) graphs $(G,o)$ with $\R$-valued labels $\{c_G(v):v\in V(G)\}\cup\{c_G(e):e\in E(G)\}$ on the vertices and edges, equipped with the topology generated by the sets 
\begin{align*}
\left\{\begin{array}{r l}
(G,o)\in\GG_\star: &\exists \varphi:B_G(o,r)\to H \textrm{ rooted (directed) graph homomorphism s.t. } \\
 &|c_G(a)-c_H(\varphi(a))|<\eps, \forall a\in V(B_G(o,r))\cup E(B_G(o,r))
\end{array}\right\},
\end{align*} 
where $\eps>0$, $r$ is any positive integer, $H$ is any finite rooted (directed) graph with labels $\{c_H(a):a\in V(H)\cup E(H)\}$ on the vertices and edges. 
A measurable function $f:\GG_\star\to \R$ is called a \emph{factor}. 

Let $G$ be a (random directed) graph, let $c:V(G)\to [0,1]$ be IID uniform random labels on the vertices and let $G(c)$ be the random labeled graph given by the labels $c$. 
The collection of random variables $\{X_a=f((G(c),a)):a\in V(G)\cup E(G)\}$ is called a \emph{factor of IID process}, if $f$ is a factor.

A random subset $M\subseteq E(G)$ is called a \emph{factor of IID (directed) matching} if there is a factor of IID process $(X_a)$ such that an edge $e$ is in $M$ if and only if $X_e=1$ and $M$ is a matching of $G$ with probability 1 with respect to the law of $G(c)$. 
\end{defi}

We note, that given a unimodular random rooted graph $(G,o)$ and a factor of IID process $(X_a)$ on $G$, the law of the labeled rooted graph $(G,(X_a),o)$ is unimodular as well. In particular, every factor of IID matching $M$ of a unimodular graph satisfies that $(G,M,o)$ is unimodular. 

\begin{lemma}\label{lem.fIIDmatching}
(1) For any locally finite graph $G$ and any $T>0$ there is a factor of IID matching $M(T)$ that has no augmenting paths of length at most $T$. 

\noindent (2) If $(G,o)$ is a random unimodular rooted graph, then $\lim_{T\to\infty}\pr_G\left(o\in V(M(T))\right)=\Em(G,o)$. 
\end{lemma}

\begin{rem}\label{rem.dirfIIDmatching}
The above lemma holds for directed graphs as well: the statements of the lemma remain true for the pre-images of the matchings $M(T)$ by the bijection defined in Remark \ref{rem.dirm}. 
\end{rem}

The proof of part 1) of Lemma \ref{lem.fIIDmatching} is similar to that of Lemma 2.2 of \cite{EL}, but for the sake of completeness we present it here. 
The main difference is that for graphs with unbounded degrees we cannot define the matchings $M(T)$ using Borel colorings, which were used in \cite{EL}. 
To handle the case of unbounded degrees we define $M(T)$ as factor of IID matchings.  
Our language is also different, although all the claims stated for Borel matchings in \cite{EL} hold for factor of IID matchings as well. 

We need the following lemma for the proof of part 2) of Lemma \ref{lem.fIIDmatching}.

\begin{lemma}\label{lem.unimodsup}
Let $(G,o)$ be a unimodular random rooted graph. Then if a unimodular matching $M$ of $G$ satisfies that there are no augmenting paths of length at most $k$, then  
\begin{align*}
\pr(o\in V(M))\geq \Em(G,o)-1/k.
\end{align*}
\end{lemma}
\begin{pf}
We show that for every $\eps$ and $k$, any unimodular matching $M$ that has no augmenting path of length at most $k$ satisfies  
\begin{align}\label{eq.meps}
\pr(o\in V(M))\geq \Em(G,o)-\eps -1/k.
\end{align}
This implies the statement of the lemma. 
Let $M_\eps$ be a fixed unimodular matching that satisfies $\Em(G,o)-\pr(o\in V(M_\eps))\leq\eps$. Consider the symmetric difference $M\bigtriangleup M_\eps$, that is a disjoint union of paths and cycles, which alternately consists of edges of $M$ and $M_\eps$  by the definition of matchings. 
We will bound $\pr(o\in V(M_\eps)\setminus V(M))$ from above by $1/k$, which implies \eqref{eq.meps} by 
\begin{align*}
\pr(o\in V(M))\geq \pr(o\in V(M_\eps))-\pr(o\in V(M_\eps)\setminus V(M)). 
\end{align*}
If a vertex $x$ of $G$ is in $V(M_\eps)\setminus V(M)$, then there is an alternating path consisting of at least $2k+2$ edges in $M\bigtriangleup M_\eps$ starting from $x$ with an edge of $M_\eps$ by the assumption on $M$. 
Define the following mass transport: let $f(x,y,(G,M\bigtriangleup M_\eps))$ be 1, if $x\in V(M_\eps)\setminus V(M)$ and $y$ is at distance at most $k-1$ from $x$ in the graph metric induced by $M_\eps\bigtriangleup M$ (there is exactly $k$ such $y$, by our previous observation on the alternating path starting from $x$). Let $f(x,y,(G,M\bigtriangleup M_\eps))$ be 0 otherwise. Note that each vertex receives mass at most 1. The labeled graph $(G,M\bigtriangleup M_\eps,o)$ is unimodular, hence we have by the Mass Transport Principle that 
\begin{align*}
k\pr(o\in V(M_\eps)\setminus V(M))&=\ex\left(\sum_{x\in V(G)}f(o,x,(G,M\bigtriangleup M_\eps))\right)\\
&=\ex\left(\sum_{x\in V(G)}f(x,o,(G,M\bigtriangleup M_\eps))\right)\leq 1.
\end{align*}
This gives the desired bound on $\pr(o\in V(M_\eps)\setminus V(M))$. 
\end{pf}

\medskip
\begin{pfo}\emph{ Lemma \ref{lem.fIIDmatching}.}
We assign to each vertex $x$ of $G$ a uniform random $[0,1]$-label
$c(x)$. First we note that with probability 1 all the labels are different, so
we can assume this property. Furthermore, we can decompose each label $c(x)$
into countably many labels $(c_{i,j}(x))_{i,j=0}^\infty$ whose joint law is IID uniform on [0,1]. 
First we construct partitions $\VV_T=\{V_{T,j}:j\geq 1\}, T\geq 1$ of $V$ such
that for each $T$ and $j$ $\inf\{\dist(x,y):x,y\in V_{T,j}\}\geq 6T$ holds. Let 
\begin{align*}
V_{T,1}:=&\left\{x\in V: c_{T,1}(x)<c_{T,1}(y) \textrm{ for every } y\in B_G(x,6T)\right\},\\
V_{T,j}:=&\left\{x\in V\setminus \left(\bigcup_{l=1}^{j-1}V_{T,l}\right): c_{T,j}(x)<c_{T,j}(y) \textrm{ for every } y\in B_G(x,6T)\right\},\quad j\geq 2. 
\end{align*}
Since the labels are uniform in $[0,1]$, we get a partition with probability
one. 

We define the matchings $M_n(T)$ in the following way. 
Let $M_0(T)=M(T-1)$ (and the empty matching if $T=1$) and let $k(n)$ be a fixed sequence that consists of positive integers and contains each of them infinitely many times.  
To define $M_n(T)$ we improve the matching $M_{n-1}(T)$ in all the balls $B(x,3T)$ with $x\in V_{T,k(n)}$: we improve using the augmenting path of length at most $T$ lying in $B(x,3T)$ with the maximal sum of $c_{T,0}$-labels of the vertices and we repeat this as long as there are short augmenting paths. The number of vertices in $B(x,3T)$ that are incident to edges of the matching increases in each step, hence we can make only a finite number of improvements in each ball. 
Since for all $n$ the balls in $\{B(x,3T): \in V_{T,k(n)}\}$ are disjoint, $M_n(T)$ is a well defined matching for every $n$ and $T$. 

Let $M(T)$ be the edge-wise limit of $M_n(T)$ as $n\to\infty$. 
We claim that $M(T)$ is well defined and has no augmenting paths of length at most $T$. 
Indeed, an edge $e=\{x,y\}$ changes its status of being in the matching or not only if there is an improvement in $B(x,3T)$. Such an improvement increase the number of vertices incident to edges of the matching in $B(x,3T)$, which is bounded above by the number of vertices in the ball, thus the number of changes is bounded above as well. 
The lack of short augmenting paths follows trivially from the construction of $M(T)$. 

We note that every factor of IID matching $M$ of a unimodular random rooted graph $(G,o)$ satisfies that $(G,M,o)$ is unimodular, hence Lemma \ref{lem.unimodsup} implies the second statement of the theorem. 
\end{pfo}

Since we do not assume the existence of a uniform bound on the degrees, we need a lemma that plays the role of Lemma 4.1 of \cite{EL}. 
\begin{lemma}\label{lem.muH}
Let $(G,o)$ be a labeled (directed) unimodular graph with law $\mu$ and finite expected degree. Then for any $\eps>0$ and any $n$ there is a $\delta$ such that if a measurable event $H$ satisfies $\mu(H)<\delta$, then $\mu(H^n)<\eps$, where $H^n:=\{(\om,x): (\om,o)\in H, \dist_\om(o,x)\leq n\}$.
\end{lemma}
\begin{pf}
Fix $\eps$ and define $D=D(\eps)$ to be the smallest positive integer that satisfies $\ex\left(\mathbf{1}_{\{\deg o>D\}}\deg o\right)<\eps/4$. 
We define the following mass transport: let $f(x,y,\om)=1$, if $(\om,x)\in H, (\om,y)\notin H, \{x,y\}\in E(\om)$ (or in the directed case $(x,y)$ or $(y,x)\in E(\om)$), and let $f(x,y,\om)=0$ otherwise. Then by the Mass Transport Principle 
\begin{align*}
\mu(H^1\setminus H)&\leq \int \sum_{x\in V(G)} f(x,o,\om) d\mu(\om,o)
= \int \sum_{x\in V(G)} f(o,x,\om) d\mu(\om,o)\\
 &\leq \ex\left( \deg o\cdot \mathbf{1}_{\{o\in H\}} \right)\\
 &\leq \ex\left( D\cdot \mathbf{1}_{\{o\in H, \deg o\leq D\}} \right)+\ex\left( \deg o\cdot \mathbf{1}_{\{o\in H, \deg o>D\}} \right)\\
 &\leq D\mu(H)+\eps/4,
\end{align*}
which is less then $\eps/2$ if $\mu(H)<\frac{\eps}{4D(\eps)}:=\eps_1$. It follows that $\mu(H^1)<\eps$. 
We define recursively $\eps_k:=\frac{\eps_{k-1}}{4D(\eps_{k-1})}$ for $k\geq 2$. Then the same argument shows that if $\mu(H)<\eps_n$, then $\mu(H^n)<\eps$.
\end{pf}

\medskip
\begin{pfo}\emph{ Theorem \ref{thm.convem}.}
First we note that by Remark \ref{rem.dirm} and 
Proposition \ref{prop.bipartite} it is enough to prove the theorem for non-directed graphs. 

Denote the law of the limit graph $(G,o)$ endowed with IID uniform labels $c(x)$ by $\mu$. 
Fix $T$ and let $\eps_T>0$ be such that if an event $H$ satisfies $\mu(H)<\eps_T$, then $\mu(H^{2T+1})<1/T$, as provided by Lemma \ref{lem.muH}. 
Let $M(T)$ be a matching as defined in Lemma \ref{lem.fIIDmatching}. 

We define the following events:  
let $\XX_0:=\{\deg_{M(T)}o=0\}$ and 
let $\XX_{i,j}$ be the event that there is an edge
$\{o,x\}\in M(T)$, such that $x$ has the
$i^{th}$ largest label among the neighbors of $o$ and $o$ has the
$j^{th}$ largest label among the neighbors of $x$. 
Note that the above events are disjoint, $\mu\left(\XX_0\cup \left(\bigcup_{i,j}\XX_{i,j}\right)\right)=1$ and if $\{x,y\}\in M(T)$ then $(G,x)\in \XX_{i,j}$ if and only if $(G,y)\in
\XX_{j,i}$. 
We can find constants $r=r(T)$ and $d=d(T)$ which satisfy the following: there are disjoint events $\YY_{i,j}, i,j\in [d]$ and $\YY_0=\left(\cup_{i,j\in[d]}\YY_{i,j}\right)^c$ determined by the labeled neighborhood of radius $r$ such that  
$\mu\left(H\right)<\eps_T$ where $H:=(\YY_0\bigtriangleup \XX_0)\cup
\left(\bigcup_{i,j\leq d}(\YY_{i,j}\bigtriangleup \XX_{i,j})\right)\cup
\left(\bigcup_{\max\{i,j\}> d}\XX_{i,j}\right)$, 
furthermore if 
$\deg_G o>d$, then $(G,o)\in \YY_0$. 
Denote by $\B(\YY_{i,j})$ the
isomorphism types of neighborhoods of radius $r$ which determine $\YY_{i,j}$. 

Now we give all vertices of $G_n$ uniform random [0,1] labels independently and denote the joint law of $G_n$ and the labels by $\mu_n$. 
We define the random matching $M_T(G_n)$ using the labels and the sets $\B(\YY_{i,j})$: 
let an edge $\{x,y\}$ be in $M_T(G_n)$ if{f} 
there is a pair $(i,j)$ such that 
$B_{G_n}(x,r)\in \B(\YY_{i,j})$, $y$ has the $j^{th}$ largest label among the neighbors of $x$, and 
$B_{G_n}(y,r)\in \B(\YY_{j,i})$, $x$ has the $i^{th}$ largest label among the neighbors of $y$. 
The edge set $M_T(G_n)$ is a matching, because the events $\B(\YY_{i,j})$ are disjoint. 
We can define a matching $M_T(G)$ of $G$ in the same way. Note, that $M_T(G)$ does not necessarily coincide with $M(T)$ but it satisfies $|\mu(o\in V(M(T)))-\mu(o\in V(M_T(G)))|<2\eps_T$ by the definition of $M_T(G)$. 
It follows by Lemma \ref{lem.unimodsup} that $\lim_{T\to\infty}\mu\big(o\in V(M_T(G))\big)=\lim_{T\to\infty}\mu\big(o\in V(M(T))\big)=\Em(G,o)$. 

Denote by $\QQ_T$ the event that there is an augmenting path for $M_T$ of length less than $T$ starting from the root. 
Let $Q_T(G_n)$ be the random set of vertices $v$ of $G_n$ such that $(G_n,v)\in\QQ_T$ 
and let $q_T(G_n):=\frac{|Q_T(G_n)|}{|V(G_n)|}$. 
The event $(G_n,x)\in Q_T$ depends on $B_{G_n}(x,r+2T+1)$ by the definition of $M_T$. 
Furthermore, in the limiting graph $G$, an augmenting path of length less than $T$ can start from $o$ only if there is a vertex $x$ on that path with $(G,x)\in H$, hence we have $\QQ_T(G,o)\subseteq H^{2T+1}$. 
It follows from the convergence $G_n\to (G,o)$ that
\begin{align*}
\lim_{n\to\infty}\ex(q_T(G_n))=\lim_{n\to\infty}\mu_n(\QQ_T(G_n,o))\leq \mu(H^{2T+1})<\frac{1}{T}, 
\end{align*}
hence $\ex(q_T(G_n))<2/T$ for $n$ large enough. 
We have by Lemma 2.1 of \cite{EL}, that 
\begin{align}\label{eq.mGn}
\frac{|M_T(G_n)|}{|V(G_n)|}\leq m(G_n)\leq
\frac{T+1}{T}\frac{|M_T(G_n)|}{|V(G_n)|}+q_T(G_n). 
\end{align}
Taking expectation in \eqref{eq.mGn} with respect to $\mu_n$, we have for $n$ large enough that 
\begin{align*}
\mu_n\left(o\in V (M_T(G_n))\right)=\ex\left(\frac{|M_T(G_n)|}{|V(G_n)|}\right)\leq \ex(m(G_n))\leq
\frac{T+1}{T}\mu_n\left(o\in V(M_T(G_n))\right)+\frac{2}{T}, 
\end{align*}
where $o$ is a uniform random vertex of $G_n$. 
Since the event $\{o\in V(M_T(G_n))\}$ depends only on the $(r(T)+1)$-neighborhood of
$x$, the convergence of the graph sequence implies $\lim_{n\to\infty}\mu_n(o\in V(M_T(G_n)))=\mu(o\in V(M_T(G)))$. 
It follows by letting $T\to \infty$ that 
$\ex(m(G_n))$ converge to $\lim_{T\to\infty}\mu(o\in M_T(G))=\Em(G,o)$. 
\end{pfo}

\subsection{Almost sure convergence of the directed matching ratio}\label{ss.asconvm}

We will examine the network models described in Subsection \ref{ss.networkmodels}. 
As referred there, each model has a local weak limit, hence Theorem \ref{thm.convem} shows that the expected values of the directed matching ratios converge. 
In this section we will show that almost sure convergence holds as well.

First we note that the local weak convergence of a sequence $G_n$ of random graphs defined on a common probability space does not imply automatically that the sequence converges almost surely in the local weak sense (see Definition \ref{def.asconvgr}), as shown by the next example. 
Let $G_n$ be the path of length $n^2$, respectively the $n\times n$ square grid, with probability 1/2-1/2. Let the joint law of the sequence $G_n$ given by the product measure. Then $G_n$ converges in the local weak sense to the infinite rooted graph $G$ which is $\Z$, respectively $\Z^2$ with probability 1/2-1/2, but there is almost surely no local weak limit of the deterministic graph sequence given by the product measure. 
If a sequence $G_n$ of finite random graphs converges almost surely in the local weak sense, then Theorem \ref{thm.convem} implies the almost sure convergence of the matching ratio, which will be the case for some of the examined sequences. 

\begin{rem}\label{rem.Skorohod}
Skorohod's Representation Theorem states that for a weakly convergent sequence $\mu_n\to \mu$ of probability measures on a complete separable metric space $S$ there is a probability space $(\Omega,\mathcal{F},\mathcal{P})$ and $S$-valued random variables $X_n$ and $X$ with laws $\mu_n$ and $\mu$ respectively, such that $X_n\to X$ almost surely. 

One could think that Skorohod's Theorem could be applied for the graph sequences that we consider, and get the convergence of the matching ratio for almost every sequence, using Theorem \ref{thm.convem}. This argument does not work for our purpose, because in Skorohod's Theorem, the coupling between the finite graphs is coming from the theorem, while in the case of the preferential attachment graphs there is given a joint probability space by construction, that contains them all. 
\end{rem}

We present two distinct methods to prove the existence of the almost sure limit of the matching ratio of a convergent graph sequence $G_n$. 
The first one can be applied for the random graph models of Section \ref{ss.networkmodels} that are defined by giving the edges independent orientations. We use this method in Subsection \ref{sss.diras} to prove part 2) of Theorem \ref{thm.convmall}. 
We show in Lemma \ref{lem.dirasconv} that if we give the edges of a converging \emph{deterministic} graph sequence uniform random orientations, then the obtained graph sequence converges almost surely in the local weak sense (see Definition \ref{def.asconvgr}) to the same limiting graph with randomly oriented edges.  
Applying this result to the sequences of Erd\H{o}s--R\'enyi random graphs and the random configuration model, which are known to converge almost surely in the non-directed case, the almost sure convergence of the matching ratio follows by Theorem \ref{thm.convem}.

We apply the approach with the second type of argument to preferential attachment graphs in Subsection \ref{sss.pa}. 
The first method does not apply for this class of graphs, because the orientations of the edges are not independent and we cannot start from an a priory almost sure convergence of the undirected graph sequence. 
We will show that the matching ratio of $G_n$  is concentrated around its expected value, which together with Theorem \ref{thm.convem} on the convergence of the mean of the matching ratio implies the almost sure convergence.

\subsubsection{Directed versions of almost surely convergent graph sequences}\label{sss.diras}

In this section we prove Part 2) of Theorem \ref{thm.convmall}. As a consequence, we have that the directed matching ratios of sequences of random regular graphs, graphs given by the random configuration model and Erd\H{o}s--R\'enyi random graphs converge almost surely, see 
Corollary \ref{cor.rrgasconv}, Theorem 
\ref{cor.randconfmasconv} and Corollary \ref{cor.ERasconv}, respectively. 

First we prove Lemma \ref{lem.dirasconv} on the almost sure convergence of a sequence of random directed graphs (see Definition \ref{def.asconvgr}) obtained from a convergent deterministic graph sequence by giving independent uniform orientation to the edges. 
This lemma implies Part 2) of Theorem \ref{thm.convmall}. 

The graph sequences examined in this section are known to converge almost surely in the undirected case. 
It follows by Part 2) of Theorem \ref{thm.convmall} that their directed matching ratios converge almost surely. 
By our Proposition 
\ref{prop.bipartite} and Theorem 2 in \cite{BLS} on the limit of the matching ratio of convergent graph sequences, one can compute the value of the limit of the directed matching ratio when the limit is a unimodular Galton--Watson tree. In Corollaries \ref{cor.rrgasconv} and \ref{cor.ERasconv} we also present the results given by this argument. 

\begin{lemma}\label{lem.dirasconv}
Let $G_n$ be a sequence of deterministic undirected graphs on $n$ vertices that converges to the random rooted graph $(G,o)$ in the local weak sense. 
Let $G^d_n$ be the sequence of random directed graphs obtained from $G_n$ by giving a random uniform orientation to each edge uniformly independently. 
Then the sequence $G^d_n$ converges almost surely in the local weak sense to $(G^d,o)$, which is the random rooted graph obtained from $(G,o)$ by orienting each edge independently. 
\end{lemma}

\begin{pfo}\emph{ Theorem \ref{thm.convmall}, Part 2).}
Consider a sequence $G_n^d$ of random directed graphs obtained by giving a uniform random orientations to the edges of a sequence of undirected random graphs $G_n$ that converges almost surely in the local weak sense to the limit graph $(G,o)$. 
We have by Lemma \ref{lem.dirasconv}, that $G_n^d$ converges almost surely in the local weak sense to the directed graph $(G^d,o)$. 
It follows by Theorem \ref{thm.convem}, that the sequence $m(G_n)$ of the matching ratios converges almost surely to $\Em(G^d,o)$.
\end{pfo}

\medskip
The proof of Lemma \ref{lem.dirasconv} essentially follows the proof of Proposition 2.2 in \cite{E09}. 
The main difference is that in that proof there were considered graphs with an uniform bound on the degrees. 

\begin{pfo}\emph{ Lemma \ref{lem.dirasconv}.}
To handle the case of unbounded degrees, we consider the following neighborhoods of the vertices: for any graph $G$ and $v\in V(G)$ denote by $B^-_{G}(v,r)$ the subgraph of $G$ obtained from $B_{G}(v,r)$ by removing all edges with both endpoint being at distance $r$ from $v$. 
Then the local weak convergence of the sequence of the finite (directed) random graphs $G_n$ to the rooted random (directed) graph $(G,o)$ is equivalent with the following: for any $r$ and any finite (directed) rooted graph $H$ we have $\lim_{n\to\infty}\pr(B^-_{G_n}(o_n,r)\simeq H)=\pr(B^-_{G}(o,r)\simeq H)$, where $o_n$ is a uniform random vertex of $G_n$. 

Fix any positive integer $r$ and any finite directed rooted graph $H^d$. 
Let $H$ be the rooted non-directed graph obtained from $H^d$ by forgetting the orientations of the edges. 
Denote by $b(G_n)$ and $b(G^d_n)$ the number of vertices $v$ of $G_n$ (respectively $G^d_n$) such that $B^-_{G_n}(v,r)\simeq H$ (resp. $B^-_{G^d_n}(v,r)\simeq H^d$). 
We show that $\pr\left(B^-_{G^d_n}(o,r)\simeq H^d\right)=\frac{b(G^d_n)}{n}$ almost surely converges to $\pr\left(B^-_{G^d}(o,r)\simeq H^d\right)$. Since this holds for any $H^d$, the lemma follows. 

Let $h$ be the probability that the graph obtained from $H$ by giving each edge a random orientation independently is isomorphic to $H^d$. 
Then $\ex(b(G^d_n))=hb(G_n)$. 
We will show that 
\begin{align}\label{eq.bconv}
\frac{b(G^d_n)}{b(G_n)}\to h \textrm{ almost surely.} 
\end{align}
The statement of the lemma follows from this, because the assumption on the convergence of $G_n$ implies that $\frac{hb(G_n)}{n}$ converges to $h\pr(B^-_G(o,r)\simeq H)=\pr(B^-_{G^d}(o,r)\simeq H^d)$. 

To show \eqref{eq.bconv}, we note that if two vertices $x$, $y$ in $G_n$ satisfy $B^-_{G_n}(x,r)\simeq B^-_{G_n}(y,r) \simeq H$ and $\dist_{G_n}(x,y)\geq 2r$, then the orientations of all the edges in $B^-_{G^d_n}(x,r)\cup B^-_{G^d_n}(y,r)$ are independent. Let $D$ be the maximum degree of the graph $H$. 
We claim that we can define a partition $(R_j^n)_{j=1}^{D^{2r}+1}$ of the set $\{x\in V(G_n): B^-_{G_n}(x,r)\simeq H\}$ such that the distance between any two points of $R_j^n$ is at least $2r$ for every $j$ and $n$. 
Indeed, if $\dist_{G_n}(x,y)$ is less than $2r$ and $B_{G_n}^-(x,r)\simeq B_{G_n}^-(y,r)\simeq H$, then there is a path of length at most $2r-1$ such that every vertex of that path has distance at most $r-1$ from the set $\{x,y\}$, and hence every vertex in the path has degree at most $D$. 
It follows, that for any fixed $x$, the number of such paths and hence the number of vertices $y$ with $\dist_{G_n}(x,y)<2r$ is at most $D^{2r}$. 
We conclude as in the proof of Proposition 2.2 in \cite{E09}. The further part of the proof is the same as the proof of that, but for the sake of completeness we present it here. 
The graph with vertex set $\{x\in V(G_n): B_{G_n}(x,r)\simeq H\}$ and edge set $\{\{x,y\}:\dist_{G_n}(x,y)<2r\}$ has maximal degree at most $D^{2r}$, thus there is a coloring of its vertices with $D^{2r}+1$ colors, that gives the partition $(R_j^n)$. 
Let $\eps>0$ and $\delta>0$ be arbitrary and let $R_1^n,\dots,R_{k(n)}^n$ be the list of the sets $R_j^n$ which satisfy $|R_j^n|\geq \eps|V(G_n)|/(D^{2r}+1)$. Denote by $b(R_j^n)$ the number of vertices $v$ in $R_j^n$ such that $B_{G_n}(v,r)\simeq H^d$. 
By the strong law of large numbers 
\begin{align*}
\left|\frac{b(R_j^n)}{|R_j^n|}-h\right|<\eps
\end{align*}
holds for all $n$ large enough and $j\leq k(n)$ with probability at least $1-\delta$, and hence we have that 
\begin{align*}
\left|\frac{b(G^d_n)}{b(G_n)}-h\right|
\leq \left|\frac{b(G^d_n)}{b(G_n)}-\frac{\sum_{j=1}^{k(n)}b(R_j^n)}{\sum_{j=1}^{k(n)}|R_j^n|}\right|
+\left|\frac{\sum_{j=1}^{k(n)}b(R_j^n)}{\sum_{j=1}^{k(n)}|R_j^n|}-h\right|
\leq \eps \frac{b(G^d_n)}{b(G_n)}+\eps\leq 2\eps 
\end{align*}
for all large enough $n$ with probability at least $1-\delta$. Since $\eps$ and $\delta$ was arbitrary, this implies \eqref{eq.bconv}. 
\end{pfo}

The directed versions of the first three graph models of Subsection \ref{ss.networkmodels} are given by orienting the edges of the non-directed versions independently. 
We use the following consequence of Theorem 3.28 in \cite{Brg} for the almost sure convergence of the directed random configuration model (see Subsection \ref{ss.networkmodels} for the definition): 
\begin{tetel}[\cite{Brg}, Theorem 3.28.]\label{thm.Bordenave}
If $G_n$ is sequence of random undirected graphs given by the random configuration model with degree distribution $\xi$ satisfying $\ex(\xi^p)<\infty$ for some $p>2$, then the sequence $G_n$ converges to $UGW(\xi)$ almost surely in the local weak sense. 
\end{tetel}

A corollary of Part 2) of Theorem \ref{thm.convmall} and Theorem \ref{thm.Bordenave} is the almost sure convergence of the sequence of graphs obtained by the random configuration model and the matching ratio of it. 
\begin{cor}[Almost sure convergence of the directed matching ratio of the random configuration model]\label{cor.randconfmasconv}
Let $G_n$ be a sequence of random directed graphs given by the random configuration
model with degree distribution $\xi$ satisfying $\ex(\xi^p)<\infty$ for some $p>2$. 
Then $G_n$ converge almost surely in the local weak sense to $UGW^d(\xi)$ and $m(G_n)$ converges almost surely to $\Em(UGW^d(\xi))$. 
\end{cor}

The sequence of random directed $d$-regular graphs is a special case of the random configuration model (with degree distribution $\xi$ being constant $d$). The connected component of the root $o'$ of the bi-partite representation $\T'_d$ has law $UGW(Binom(d,1/2))$, hence we have the following: 
\begin{cor}[Almost sure convergence of the directed matching ratios of directed random regular graphs]\label{cor.rrgasconv}
Let $G_n$ be the sequence of random $d$-regular graphs on $n$ vertices with randomly oriented edges. Then the matching ratios converge almost surely to the constant 
\begin{align*}
\lim_{n\to\infty}m(G_n)=\Em\big(UGW(Binom(d,1/2))\big).
\end{align*}
\end{cor}

For directed Erd\H{o}s--R\'enyi graphs one can compute the exact value of the almost sure limit of the matching ratio, using the results of \cite{KS} or Theorem 2 in \cite{BLS}. 

\begin{cor}\label{cor.ERasconv}
Let $\mathcal{G}_{n,2c/n}^d$ be a sequence of directed Erd\H{o}s--R\'enyi graphs with parameter $2c$. Then almost surely 
\begin{align}\label{eq.tc}
\lim_{n\to\infty}m(\mathcal{G}_{n,2c/n})=1-\frac{t_{c}+e^{-ct_{c}}+ct_{c}e^{-ct_{c}}}{2} 
\end{align}
where $t_{c}\in(0,1)$ is the smallest root of $t=e^{-ce^{-ct}}$.
\end{cor}
\begin{pf}
According to Subsection \ref{ss.networkmodels} and Lemma \ref{lem.dirasconv}, the sequence of directed Erd\H{o}s--R\'enyi random graphs converge almost surely in the local weak sense to $UGW^d(Poisson(2c))$, and hence $\lim_{n\to\infty}m(\mathcal{G}^d_{n,2c/n})=\Em(UGW^d(Poisson(2c)))$.  
The connected component of the root in the bipartite representation of $UGW^d(Poisson(2c))$ has law $UGW(Poisson(c))$, which is the almost sure local weak limit of the non-directed Erd\H{o}s--R\'enyi random graphs $\mathcal{G}_{n,c/n}$ with parameter $c$. It is known (see \cite{KS} or Theorem 2 in \cite{BLS}), that for this graph sequence $\lim_{n\to\infty}m(\mathcal{G}_{n,c/n})$ equals the right hand side of \eqref{eq.tc} almost surely. 
By Remark \ref{rem.infinitem} we have $\lim_{n\to\infty}m(\mathcal{G}_{n,c/n})=\Em(UGW(Poisson(c)))=\Em(UGW^d(Poisson(2c)))$. This proves \eqref{eq.tc}. 
\end{pf}

\subsubsection{Preferential attachment graphs}\label{sss.pa}

In this section we show that the directed matching ratio of a graph sequence given by the preferential attachment rule converges almost surely, see Theorem \ref{thm.asconvm}. 
The orientations of the edges of this class of graphs are given naturally by the recursive definition, and differ significantly from the independent random orientation. Thus we cannot apply the results of Section \ref{sss.diras}. 
This sequence also does not satisfy the assumption of \cite{LSB} that the distributions of the in- and out-degrees are the same (which was assumed to simplify the calculations made there), hence the value of the limit of the matching ratio for this class was not examined in that paper. 
However, the almost sure convergence of the directed matching ratios holds for this class of graph sequences as well, as we show in the next theorem. 

\begin{tetel}\label{thm.asconvm}
Let $G_n$ be a random graph sequence obtained by the preferential attachment rule. 
Then $\lim_{n\to\infty}m(G_n)=\lim_{n\to\infty}\ex(m(G_n))$ almost surely. 
\end{tetel}

We will prove the concentration of the matching ratios around their expected value, which together with the results of \cite{BBCS} on the local weak convergence of $G_n$ and Theorem \ref{thm.convem} on the convergence of the mean of the matching ratio implies the statement. 

\begin{rem}
It follows from the concentration shown in the proof, that the almost sure local weak convergence holds for any joint law of the graphs $G_n$. 
\end{rem}

\begin{pfo}\emph{ Theorem \ref{thm.asconvm}.}
Fix $n$ and denote by $G_n(k)$ the subgraph of $G_n$ spanned by the vertices
$\{1,\dots,k\}$. Let 
\begin{align}
X_k:=\ex\left(|M_{max}(G_n)|\Big|G_n(k)\right)-\ex\left(|M_{max}(G_n)|\Big|G_n(k-1)\right).
\end{align}
We will show that $|X_k|\leq 2r$ almost surely for all $k\in[n]$. 
Since $Y_k:=\ex\left(|M_{max}(G_n)|\Big|G_n(k)\right)$ is a martingale,
we can apply the Azuma--Hoeffding inequality (Theorem \ref{thm.AH}) to the random variables $X_k$. 
It follows that for any $c>0$ we have 
\begin{align*}
\pr\left( |m(G_n)-\ex(m(G_n))|>c \right)& 
=\pr\left( |X_1+\dots +X_n|>cn \right)\\
&\leq 2\exp\left\{ -\frac{(cn)^2}{2\sum_{k=1}^n \|X_k\|_\infty^2} \right\}\\
&\leq 2\exp\left\{ -\frac{c^2n^2}{8nr^2} \right\}. 
\end{align*}
Since $\lim_{k\to\infty}\ex(m(G_k))$ exists by Theorem \ref{thm.convem}, for $n$ large enough to satisfy $|\ex(m(G_n))-\lim_{k\to\infty}\ex(m(G_k))|<c/2$ we have 
\begin{align*}
\pr\left( |m(G_n)-\lim_{k\to\infty}\ex(m(G_k))|>c \right)& 
\leq\pr\left( |m(G_n)-\ex(m(G_n))|>\frac{c}{2} \right)\\
&\leq 2\exp\left\{ -\frac{c^2n}{32r^2} \right\}. 
\end{align*}
It follows that these probabilities are summable in $n$ for every $c>0$, which implies the almost sure convergence of $m(G_n)$ by the Borel-Cantelli lemma.

What remains to show is that for any fixed pair of
directed graphs $F$ and $F'$ on the vertex set $[k]$ with $F(k-1)=F'(k-1)$, the inequality 
\begin{align}
\left|\ex\big(|M_{max}(G_n)|\big| G_n(k)=F \big.\big)-\ex\left(|M_{max}(G_n)|\left| G_n(k)=F' \right.\right)\right|\leq 2r
\label{eq.Mest}
\end{align} 
holds. This implies $|X_k|\leq 2r$. 

Fix $F$ and $F'$ as above. For any possible configuration of $G_n$, denote by 
\begin{align*}
h(G_n):=\left\{(\ell,j)\in E(G_n): \ell>k, (k,j)\notin E(F)\cup E(F')\right\}
\end{align*}
the subset of the edges of $G_n$ with tails in $\{k+1,\dots,n\}$ that do not
have a common head with the edges in the graphs $F$ or $F'$ with tail $k$.
The proof of inequality \eqref{eq.Mest} is based on two observations:
first, by the definition of the preferential attachment graph, 
the distribution of $h(G_n)$ conditioned on $\{G_n(k)=F\}$ is
the same as conditioned on $\{G_n(k)=F'\}$ (note the symmetry in $F$ and $F'$ in the definition of $h(G_n)$). 
Second, for any configuration of $G_n$ with $G_n(k)=F$, the size of the maximal matching changes by at most
$2r$ if we 
fix $h(G_n)$, set $G_n(k):=F'$ 
and vary arbitrary the heads of the edges with tails in $\{k+1,\dots ,n\}$ that are not in $h(G_n)$. 
This follows from Lemma \ref{lem.adde} by the following argument. 
For any fixed $H$ 
we obtain any graph in the set $\{G_n:G_n(k)=F,h(G_n)=H\}$ by adding 
new edges with heads in the set $\{j: (k,j)\in E(F)\cup E(F')\}$ of size at
most $2r$ 
to the graph $G_H$ with
$V(G_H):=[n]$ and $E(G_H):=E(F(k-1))\cup H$. 
It follows from Lemma \ref{lem.adde} that 
\begin{align}\label{eq.Mest2}
|M_{max}(G_H)|
\leq \ex\big(|M_{max}(G_n)|\big|G_n(k)=F,h(G_n)=H\big)
\leq |M_{max}(G_H)|+2r,
\end{align}
and the same holds with $F'$ in the place of $F$. This proves the second observation. 

Using the first observation and \eqref{eq.Mest2} the left hand side of \eqref{eq.Mest} can be estimated from above by 
\begin{align}
\sum_{H} &\left|
\ex\left(|M_{max}(G_n)|\Big|G_n(k)=F,h(G_n)=H\right)\pr\left(h(G_n)=H\Big|G_n(k)=F\right)\right.\nonumber\\
&\left. - \ex\left(|M_{max}(G_n)|\Big|G_n(k)=F',h(G_n)=H\right)\pr\left(h(G_n)=H\Big|G_n(k)=F'\right) \right|\nonumber \\
\leq &\sum_{H} \pr\left(h(G_n)=H\Big|G_n(k-1)=F(k-1)\right)\cdot\nonumber\\
& \cdot\left| \ex\left(|M_{max}(G_n)|\Big|G_n(k)=F,h(G_n)=H\right)-
\ex\left(|M_{max}(G_n)|\Big|G_n(k)=F',h(G_n)=H\right) \right| \nonumber\\
\leq &\sum_{H} \pr\left(h(G_n)=H\Big|G_n(k-1)=F(k-1)\right)\cdot2r \nonumber\\
=&2r\nonumber
\end{align}
\end{pfo}

\end{document}